\documentclass[halfparskip,pdftex]{scrartcl}

\usepackage{amsmath,amsthm,amssymb}
\usepackage{ifpdf}

\usepackage[round]{natbib} 

\usepackage{calc, graphicx} 

\usepackage{bbm}
\newcommand{\ind}{\ensuremath{\mathbbm{1}}}

\usepackage{algorithmicx,algpseudocode}
\algnewcommand\algorithmicto{\textbf{to}}
\algrenewtext{For}[3]%
  {\algorithmicfor\ #1 $\gets$ #2 \algorithmicto\ #3 \algorithmicdo}

\newtheorem{sz}{Theorem}[section]
\newtheorem{lm}[sz]{Lemma}
\newtheorem{cor}[sz]{Corollary}
{\theoremstyle{definition}
  \newtheorem{bem}[sz]{Remark}
}

\newcommand{\Nat}{\ensuremath{\mathbb{N}}}
\newcommand{\Real}{\ensuremath{\mathbb{R}}}
\newcommand{\Zs}{\ensuremath{\mathbb{Z}}}
\newcommand{\Erw}[1]{\ensuremath{\mathbb{E}\!\left[#1\right]}}
\newcommand{\Ew}{\ensuremath{\mathbb{E}}}
\newcommand{\Ws}[2][]{\ensuremath{\mathbb{P}_{#1}\!\left[#2\right]}}
\newcommand{\Pb}{\ensuremath{\mathbb{P}}}

\newcommand{\abs}[1]{\ensuremath{\left\lvert#1\right\rvert}}
\newcommand{\norm}[2][]{\ensuremath{\left\lVert#2\right\rVert}_{#1}}

\newcommand{\coleq}{\ensuremath{\mathop{:}\!\!=}}
\newcommand{\eqcol}{\ensuremath{=\!\!\mathop{:}}}
\newcommand{\minl}[1][p]{\ensuremath{\ell_{#1}}}
\newcommand{\kol}{\ensuremath{\rho}}
\newcommand{\eqd}{\ensuremath{\stackrel{d}{=}}}
\newcommand{\dist}{\ensuremath{\mathcal{L}}}

\newcommand{\An}[1][n]{A^{(#1)}}
\newcommand{\bn}[1][n]{b^{(#1)}}
\newcommand{\Xn}[1][n]{X_{#1}}
\newcommand{\Xtn}[1][n]{\tilde X_{#1}}
\newcommand{\RAn}[1][n]{R_{A}(#1)}
\newcommand{\Rbn}[1][n]{R_{b}(#1)}
\newcommand{\RXn}[1][n]{R_{X}(#1)}
\newcommand{\Rlpn}[1][n]{R(#1)}
\newcommand{\ximax}[1][p]{\mbox{$\xi\hspace{0.2ex}$}_{\!#1}}
\newcommand{\CA}[1][p]{C_{A}}
\newcommand{\Cb}[1][p]{C_{b}}
\newcommand{\CX}[1][p]{C_{X}}

\newcommand{\f}[1][x]{f_{X}(#1)}
\newcommand{\ff}{f_{X}}
\newcommand{\Dbeta}[2]{\ensuremath{\mathrm{beta}\!\left(#1,#2\right)}}
\newcommand{\DBern}[1]{\mathrm{Bernoulli\!\left(#1\right)}}
\newcommand{\flr}[1]{\left\lfloor#1\right\rfloor}
\newcommand{\ceil}[1]{\left\lceil#1\right\rceil}
\newcommand{\deltaf}[1]{\ensuremath{J_{f_{#1}}}}
\newcommand{\modcont}{\ensuremath{\Delta}}
\newcommand{\cpp}{C\hspace{-.05em}\raisebox{.4ex}{\tiny\bf +}\hspace{-.05em}\raisebox{.4ex}{\tiny\bf +}}

\newcommand{\discrete}[2][\!n]{\left\langle#2\right\rangle_{#1}}
\newcommand{\discretet}[2][\!n]{\bigl\langle#2\bigr\rangle_{#1}}

\let\theta\vartheta

\let\phi\varphi

\let\rho\varrho

\let\epsilon\varepsilon

\let\tilde\widetilde
\let\ln\log

\title{Approximating Perpetuities} 
\author{Margarete Knape\thanks{Email:
    knape@math.uni-frankfurt.de} \ and Ralph Neininger\thanks{Supported by an Emmy Noether Fellowship of the Deutsche
   Forschungsgemeinschaft, \newline Email: neiningr@math.uni-frankfurt.de}\\
Department~for Mathematics and Computer Science\\ J.W.~Goethe-University
Frankfurt~a.M.\\60054 Frankfurt a.M.\\
Germany}
\begin{document}

\maketitle

\begin{abstract}
  We propose and analyze an algorithm to approximate distribution
  functions and densities of perpetuities. Our algorithm refines an
  earlier approach based on iterating discretized versions of the fixed
  point equation that defines the perpetuity. We significantly reduce
  the complexity of the earlier algorithm. Also one particular
  perpetuity arising in the analysis of the selection algorithm
  Quickselect is studied in more detail. Our approach works well for
  distribution functions.  For densities we have weaker error bounds
  although computer experiments indicate that densities can also be
  approximated well.
\end{abstract}


\section{Introduction}

A perpetuity is a random variable $X$ in $\Real$
that satisfies the stochastic fixed-point equation
\begin{equation}
  \label{eq:fixedpoint}
  X\eqd AX+b,
\end{equation}
where the symbol $\eqd$ denotes that left and right
hand side in \eqref{eq:fixedpoint} are identically
distributed and where $(A,b)$  is a vector of
random variables being independent of $X$, whereas
dependence between $A$ and $b$ is allowed.

Perpetuities arise in various different contexts: In discrete
mathematics, perpetuities come up as the limit distributions of certain
count statistics of decomposable combinatorial structures such as random
permutations or random integers. In these areas, perpetuities often
arise via relationships to the GEM and Poisson-Dirichlet distributions;
see \citet*{ArBaTa03} for perpetuities, GEM and Poisson-Dirichlet
distribution in the context of combinatorial structures; see
\citet{DonGrim93} for occurrences in probabilistic number theory.  In
the probabilistic analysis of algorithms, perpetuities arise as limit
distributions of certain cost measures of recursive algorithms such as
the selection algorithm Quickselect, see e.g.~\citet{HwTsai02} or
\citet*{MaMoSmy95}.  In insurance and financial mathematics, a perpetuity
represents the value of a commitment to make regular payments, where $b$
represents the payment and $A$ a discount factor both being subject to
random fluctuation; see, e.g.~\citet{GolMal00} or \citet*[Section
8.4]{EmKlMi97}.

As perpetuities are given implicitly by their fixed-point
characterization \eqref{eq:fixedpoint}, properties of their
distributions are not directly amenable. Nevertheless, various questions
about perpetuities have already been settled.  Necessary and sufficient
conditions on $(A,b)$ for the fixed-point equation \eqref{eq:fixedpoint}
to uniquely determine a probability distribution are discussed in
\citet{Ver79} and \citet{GolMal00}. The types of distributions possible
for perpetuities have been identified in \citet*{AIR07}.  Tail behavior
of perpetuities has been studied for certain cases in \citet{GolGru96}.

In the present article, we are interested in the central
region of the distributions. The aim is to algorithmically approximate
perpetuities, in particular their distribution functions and their
Lebesgue densities (if they exist).

For this, we apply and refine a method proposed in \citet{DevNein02}
that was originally designed for random variables $X$ satisfying
distributional fixed-point equations of the form
\begin{equation}\label{glnn1}
X\eqd \sum_{r=1}^K A_r X^{(r)}+b,
\end{equation}
where $X^{(1)},\ldots,X^{(K)},(A_1,\ldots,A_K,b)$ are independent with
$X^{(r)}$ being identically distributed as $X$ for $r=1,\ldots, K$ and
random coefficients $A_1,\ldots,A_K$, $b$, and $K\ge 2$.

The case of perpetuities, i.e., $K=1$, structurally differs from the
cases $K\ge 2$: The presence of more than one independent copy of $X$ on
the right hand side in \eqref{glnn1} often has a smoothing effect so
that under mild additional assumptions on $(A_1,\ldots,A_K,b)$ the
existence of smooth Lebesgue densities of $X$ follows, see
\citet{FilJans00} and \citet{DevNein02}. On the other hand, the case
$K=1$ often leads to distributions ${\cal L}(X)$ that have no smooth
Lebesgue density; an example is discussed in Section \ref{sec:keyEx}.

Our basic approach to approximate perpetuities is as follows:
A random variable $X$ satisfies the distributional identity
\eqref{eq:fixedpoint} if and only if its distribution is a fixed-point
of the map $T$ on the space $\mathcal{M}$ of probability distributions,
given by
\begin{equation}
  T: \mathcal{M}\rightarrow\mathcal{M},\ \mu\mapsto\dist(AY+b),
\end{equation}
where $Y$ is independent of $(A,b)$, and $\dist(Y)=\mu$. Under the
conditions $\|A\|_p<1$ and $\|b\|_p<\infty$ for some $p\ge 1$, which we
assume throughout the paper, this map is a contraction on certain
complete metric subspaces of $\mathcal{M}$. Hence, $\dist(X)$ can be
obtained as limit of iterations of $T$, starting with some distribution
$\mu_{0}$.

However, it is not generally possible to algorithmically compute the
iterations of $T$ exactly.  We therefore use discrete approximations
$(A^{(n)},b^{(n)})$ of $(A,b)$, which become more accurate for
increasing $n$, to approximate $T$ by a mapping $\tilde T^{(n)}$,
defined by
\begin{equation*}
  \tilde T^{(n)}: \mathcal{M}\rightarrow\mathcal{M},\
  \mu\mapsto\dist\!\left(A^{(n)}Y+b^{(n)}\right),
\end{equation*}
where again $Y$ is independent of
$(A^{(n)},b^{(n)})$ and $\dist(Y)=\mu$.

To allow for an efficient computation of the approximation, we impose
a further discretisation step $\discrete{\cdot}$, introduced in
Section~\ref{sec:convRates}, defining
\begin{equation*}
  T^{(n)}:\mathcal{M}\rightarrow\mathcal{M},\
  \mu\mapsto\dist\!\left(\discrete{A^{(n)}Y+b^{(n)}}\right),
\end{equation*}
where $Y$ is independent of $(A^{(n)},b^{(n)})$ and $\dist(Y)=\mu$.

In Section~\ref{sec:convRates}, we give conditions
for $T^{(n)}\circ T^{(n-1)}\circ\dots\circ
T^{(1)}(\mu_{0})$ to converge to the perpetuity given as the solution of
\eqref{eq:fixedpoint}. To this aim, we derive a
rate of convergence in the minimal $L_{p}$ metric
$\minl[p]$, defined on the space $\mathcal{M}_{p}$
of probability measures on $\Real$ with finite absolute
$p$th moment by
\begin{equation}\label{eq:defminl}
  \minl(\nu,\mu)\coleq
  \inf\left\{
    \norm[p]{V-W\bigr.}:\ \dist(V)=\nu,\dist(W)=\mu
  \right\}, \quad\text{ for }\nu,\mu\in\mathcal{M}_{p},
\end{equation}
where $\norm[p]{\cdot}$ denotes the $L_{p}$-norm of random variables.  To
get an explicit error bound for the distribution function, we then
convert this into a rate of convergence in the Kolmogorov metric $\kol$,
defined by
\begin{equation*}
  \kol(\nu,\mu)\coleq
   \sup_{x\in\Real}\abs{F_{\nu}(x)-F_{\mu}(x)\big.},
\end{equation*}
where $F_{\nu},F_{\mu}$ denote the distribution functions of
$\nu,\mu\in\mathcal{M}_{p}$.  This implies explicit rates of convergence for
distribution function and density, depending on the corresponding moduli
of continuity of the fixed-point.

For these moduli of continuity we find global bounds for perpetuities
with \mbox{$b\equiv1$} in Section~\ref{sec:ax1}. For cases with random
$b$, we have to derive these moduli of continuity individually. One
example, connected to the selection algorithm Quickselect, is worked out
in detail in Section~\ref{sec:keyEx}.

We analyze the complexity of our approach in Section \ref{sec:complex}.
As a measure for the complexity of the approximations for distribution
function and density, we use the number of steps needed to obtain an
approximation that has distance, in supremum norm, of at most $1/n$ to
the true function. Although we generally follow the approach in
\citet{DevNein02}, we can improve the complexity significantly by using
different discretisations.  For the approximation of the distribution
function to an accuracy of $1/n$ in a typical case, we obtain a
complexity of $O(n^{1+\epsilon})$ for any $\epsilon>0$. In comparison,
the algorithm described in \citet{DevNein02}, which originally was
designed for fixed-point equations of type \eqref{glnn1} with $K\ge 2$,
would lead to a complexity of $O(n^{4+\epsilon})$, if applied to our
cases. For the approximation of the density to an accuracy of $1/n$, we
obtain a complexity of $O(n^{1+1/\alpha+\epsilon})$ for any $\epsilon>0$
in the case of $\alpha$-H\"older continuous densities, cf.
Corollary~\ref{cor:complex}.

An extended abstract of this article appeared in \citet{KnN07}.

\section{Discrete approximation and convergence}
\label{sec:convRates}

Recall that our basic assumption in equation \eqref{eq:fixedpoint} is
that $\norm[p]{A}<1$ and $\norm[p]{b}<\infty$ for some $p\geq1$.
To obtain an algorithmically computable approximation of the solution of
the fixed-point equation \eqref{eq:fixedpoint}, we use an approximation
of the sequence defined as follows: We replace $(A,b)$ by a sequence of
independent discrete approximations $(\An,\bn)$, converging to $(A,b)$
in $p$th mean for $n\rightarrow\infty$. To reduce the complexity, we
introduce a further discretisation step $\discrete{\cdot}$, which
reduces the number of values attained by $X_{n}$:
 \begin{align}
   \label{eq:defXtn}
   \Xn[0]\coleq\discretet[0]{\Ew X},\quad
   \Xtn\coleq \An\Xn[n-1]+\bn,\quad
   \Xn\coleq\discretet{\Xtn},\quad
   n\geq1.
 \end{align}
 We assume that the discretisations $\An$, $\bn$ and $\discrete{\cdot}$
 satisfy
\begin{align}
  \label{eq:approxR}
  \norm[p]{\An-A}\leq \RAn, \quad
  \norm[p]{\bn-b}\leq \Rbn, \quad
  \norm[p]{\discretet{\Xtn}-\Xtn}\leq \RXn,
\end{align}
for some error functions $R_{A}$, $R_{b}$ and $R_{X}$, which we specify
later.

Furthermore, we assume that there exists some $\ximax<1$, such that for
all $n\ge 1$,
\begin{equation}
  \label{eq:xiA}
  \norm[p]{\An}\!\leq \ximax,
\end{equation}
which in applications is easy to obtain, since $\|A\|_p<1$.

By arguments similar to those used in \citet{FilJans02} and
\citet{DevNein02} we obtain the following convergence rates for the
approximations $\Xn$ to converge to the corresponding characteristics of the
fixed-point $X$. We use the shorthand notation
$\minl[p](X,Y)\coleq\minl[p](\dist(X),\dist(Y))$.

\begin{lm}\label{lm:lpR}
  Let $\left(\Xn\right)_{n\in\Nat_{0}}$ be defined by \eqref{eq:defXtn} and
  $\ximax$ as in \eqref{eq:xiA}.  Then
  \begin{equation}\label{eq:lpR}
    \minl(\Xn,X)
      \leq \ximax^{n}\,\norm[p]{X-X_{0}}
      + \sum_{i=0}^{n-1}\ximax^{i}\,\Rlpn[n-i],
  \end{equation}
  where $\Rlpn\coleq\RXn + \RAn\norm[p]{X} + \Rbn$ for the error
  functions in \eqref{eq:approxR}.
  \begin{proof}
    We have 
    \begin{align}\label{eq:minlXnX}
      \minl(\Xn,X)&\leq \minl(\Xn,\Xtn)+\minl(\Xtn,X)\nonumber\\
      &\leq \norm[p]{\discretet{\Xtn}-\Xtn}+\minl(\Xtn,X).
    \end{align}
    The first summand is bounded by \eqref{eq:approxR} and for the
    second summand we have
    \begin{align*}
      \minl(\Xtn,X)  
      &\leq \norm[p]{\Xtn-X} =\norm[p]{\An \Xn[n-1]+\bn- AX-b}\\
      &\leq \norm[p]{\An\Xn[n-1]-AX} + \norm[p]{\bn-b}\\
      &= \norm[p]{\An(\Xn[n-1]-X) - (A-\An)X} + \norm[p]{\bn-b}\\
      &\leq \norm[p]{\An}\norm[p]{\Xn[n-1]-X} + 
            \norm[p]{A-\An}\norm[p]{X} +
            \norm[p]{\bn-b},
    \end{align*}
    where in the last step we use that $\An$ and $(\Xn[n-1]-X)$ as
    well as $(A-\An)$ and $X$ are independent by assumption.

    Now we use that the infimum in the definition of $\ell_{p}$ in
    \eqref{eq:defminl} is attained and assume additionally, that
    $\Xn[n-1]$ and $X$ are chosen with
    $\norm[p]{\Xn[n-1]-X}=\minl(\Xn[n-1],X)$.  Combining this with
    \eqref{eq:minlXnX} and using the bounds given in
    \eqref{eq:approxR} and \eqref{eq:xiA}, we obtain
    \begin{equation*}
      \minl(\Xn,X)\leq \RXn + 
      \ximax\;\minl(\Xn[n-1],X) + \RAn\norm[p]{X} + \Rbn,
    \end{equation*}
    and the claim then follows by induction.
  \end{proof}
\end{lm}

To make these estimates explicit we have to specify bounds for
$\RAn$, $\Rbn$, and $\RXn$. We do so in two different ways, one
representing a polynomial discretisation of the corresponding random
variables and one representing an exponential discretisation. Better
asymptotic results are obtained by the latter one.

\begin{cor}
  \label{cor:lpR}
  Let $\Xn, n\in\Nat_{0}$ be defined by \eqref{eq:defXtn} and
   $\ximax$ as in \eqref{eq:xiA}, and assume
  \begin{equation*}
    \RAn\leq \CA\,\dfrac{1}{n^{r}},\qquad
    \Rbn\leq \Cb\,\dfrac{1}{n^{r}}, \qquad
    \RXn\leq \CX\,\dfrac{1}{n^{r}},
  \end{equation*}
  for some $r\geq1$. Then, we have
  \begin{equation*}
    \minl(\Xn,X)\leq C_{r}\,\frac{1}{n^{r}},
  \end{equation*}
  where
  \begin{equation}
    \label{eq:defCr}
    C_{r}:=
      \frac{r^{r}\,\norm[p]{X-\Xn[0]}}{\Bigl(e\,\ln\!\left(1/\ximax\right)\Bigr)^{r}} +
      \frac{r!\left(\CX + \Cb + \CA\,\norm[p]{X}\right)}%
           {\left(1-\ximax\right)^{r+1}}.
  \end{equation}
  \begin{proof}
    Using Lemma~\ref{lm:lpR} we get
    \begin{equation}\label{eq:lpr}
      \minl(\Xn,X)
      \leq \ximax^{n}\norm[p]{X-\Xn[0]}+
      (\CX+\CA\norm[p]{X}+\Cb)
      \sum_{i=0}^{n-1}\frac{\ximax^{i}}{\left(n-i\right)^{r}}.
    \end{equation}
    For the first summand, we use that the function $x\mapsto
    x^{r}\ximax^{x}$ has its maximum at $x=r/\ln(1/\ximax)$.

    To see that the second summand is of order $n^{-r}$, note that
    $1/(n-i)\leq (i+1)/n$ for all $n\geq1$ and $0\leq i\leq n-1$.
    This implies that for $\ximax<1$,
    \begin{align*}
      \sum_{i=0}^{n-1}\frac{\ximax^{i}}{\left(n-i\right)^{r}} &\leq
      \frac{1}{n^{r}}\sum_{i=0}^{n-1}\left(i+1\right)^{r}\ximax^{i}\\&
      \leq 
      \frac{1}{n^{r}}\sum_{i=0}^{\infty}(i+r)(i+r-1)\dotsm(i+1)\ximax^{i}\\
      &= \frac{r!}{\left(1-\ximax\right)^{r+1}}\,\frac1{n^{r}},
    \end{align*}
    where the last equality is obtained by differentiating the geometric
    series $r$ times.
  \end{proof}
\end{cor}

\begin{bem}
  In Corollary \ref{cor:lpR}, we are merely interested in the order of
  magnitude of $\ell_p(\Xn,X)$ without a sharp estimate of the
  constant $C_r$. When evaluating the error in an explicit example, we
  can evaluate \eqref{eq:lpr} directly to obtain sharper estimates.
\end{bem}

\begin{cor}
  \label{cor:lpRgamman}
  Let $\Xn, n\in\Nat_{0}$ be defined by \eqref{eq:defXtn} and $\ximax$ as in
  \eqref{eq:xiA}, and assume
  \begin{equation*}
    \RAn\leq \CA\,\dfrac{1}{\gamma^{n}},\qquad
    \Rbn\leq \Cb\,\dfrac{1}{\gamma^{n}}, \qquad
    \RXn\leq \CX\,\dfrac{1}{\gamma^{n}},
  \end{equation*}
  for some $1<\gamma<1/\ximax$. Then, we have
  \begin{equation*}
    \minl(\Xn,X)\leq C_{\gamma}\,\frac{1}{\gamma^{n}},
  \end{equation*}
  where
  \begin{equation}
    \label{eq:defCgamma}
    C_{\gamma}:=
      \norm[p]{X-\Xn[0]} +
      \frac{\left(\CX + \Cb + \CA\,\norm[p]{X}\right)}{1-\ximax\gamma}.
  \end{equation}
  \begin{proof}
    Using Lemma~\ref{lm:lpR} we get 
    \begin{equation}\label{eq:lprgamman}
      \minl(\Xn,X)
      \leq \ximax^{n}\norm[p]{X-\Xn[0]}+
      (\CX+\CA\norm[p]{X}+\Cb)\gamma^{-n}
      \sum_{i=0}^{n-1}\ximax^{i}\gamma^{i},
    \end{equation}
    and the assumption on $\gamma$ implies that both summands are
    $O(\gamma^{-n})$ with the constant given in the lemma.
  \end{proof}
\end{cor}

\begin{lm}
  \label{lm:kolBound}
  Let $\Xn$ and $C_{r}$ be as in Corollary~\ref{cor:lpR}
  and $X$ have a bounded density $\ff$. Then, the distance in the
  Kolmogorov metric can be bounded by
  \begin{equation}\label{eq:kolBoundPoly}
    \kol(\Xn,X)\leq
    \left(
      C_{r}\left(p+1\right)^{1/p}\norm[\infty]{\ff}
    \right)^{p/(p+1)}
     n^{-rp/(p+1)}.
  \end{equation}
Similarly, for $\Xn$ and $C_{\gamma}$ as in
Corollary~\ref{cor:lpRgamman}, we have
  \begin{equation}\label{eq:kolBoundExp}
    \kol(\Xn,X)\leq
    \left(
      C_{r}\left(p+1\right)^{1/p}\norm[\infty]{\ff}
    \right)^{p/(p+1)}
    \gamma^{pn/(p+1)}.
  \end{equation}
  \begin{proof}
    We use Lemma~5.1 in \citet{FilJans02}, which states, that for $X$
    with bounded density $f_{X}$ and any $Y$,
    \begin{equation*}
      \kol(Y,X)\leq\left(\left(p+1\right)^{1/p}\norm[\infty]{f_{X}}
        \,\minl(Y,X)\right)^{p/(p+1)}\quad\text{for }p\geq1. 
    \end{equation*}
    Using Corollaries~\ref{cor:lpR} and \ref{cor:lpRgamman}
    respectively, we get the stated result.
  \end{proof}
\end{lm}

\begin{bem}
  In some cases, we can give a similar bound, although the density of
  $X$ is not bounded or no explicit bound is known. Instead, it is
  sufficient to have a bound for the modulus of continuity of the
  distribution function $F_{X}$ of $X$, cf.~\citet{Kn06}.
\end{bem}

To approximate the density of the fixed-point, we define
\begin{equation}
  \label{eq:approxfn}
  f_{n}(x)=\frac{
    F_{n}(x+\delta_{n})-F_{n}(x-\delta_{n})
  }{2\delta_{n}},
\end{equation}
where $F_{n}$ is the distribution function of $\Xn$. For this
approximation we can give a rate of convergence, depending on the
modulus of continuity of the density of the fixed-point, which is
defined by
\begin{equation*}
  \modcont_{\ff}(\delta)\coleq\!\!
  \sup_{\substack{u,v\in\Real\\ \abs{u-v}\leq\delta}}\!
    \abs{\f[u]-\f[v]\bigr.},\quad \delta\geq0.
\end{equation*}

\begin{lm}\label{lm:difffnf}
  Let $X$ have a density $\ff$ and let $\Xn, n\in\Nat_{0}$ be defined by
  \eqref{eq:defXtn}.  Then, for $f_{n}$ defined by \eqref{eq:approxfn}
  and all $\delta_n>0$,
  \begin{equation*}
    \norm[\infty]{f_{n}-\ff\bigl.}\leq
    \frac{1}{\delta_n}\;\kol(\Xn,X)
    + \modcont_{\ff}\!\left(\delta_{n}\right).
  \end{equation*}
  \begin{proof} 
    For any $x$, we have
    \begin{align*}
      \abs{f_{n}(x)-\ff(x)\bigl.}  &\leq \abs{\frac
        {F_{n}(x+\delta_n)-F_{n}(x-\delta_n)} {2\delta_n}-\frac
        {F(x+\delta_n)-F(x-\delta_n)}
        {2\delta_n}}+\\[1ex]
      &\qquad\qquad + \abs{\frac{F(x+\delta_n)-F(x-\delta_n)}
        {2\delta_n}-\ff(x)}   \\[2ex]
      &\leq \frac{1}{\delta_n}\;\kol(\Xn,X) +\frac{1}{2\delta_n}
      \int_{-\delta_n}^{\delta_n}\!\abs{\ff(x+y)-\ff(x)\big.}dy\\
      &\leq\frac{1}{\delta_{n}}\kol(\Xn,X) +\frac1{\delta_{n}}
      \int_{0}^{\delta_{n}}\!\!\!\modcont_{f_{X}}\!(y)\,dy.
    \end{align*}
    The assertion follows since $\modcont_{\ff}$ is monotonically increasing.
  \end{proof}
\end{lm}

\begin{cor}\label{cor:boundApproxDens}
  Let $X$ have a bounded density $\ff$, which is H\"older continuous
  with exponent $\alpha\in(0,1]$. For polynomial discretisation $X_{n}$
  and $C_{r}$ as in Corollary~\ref{cor:lpR} and $f_{n}$ defined by
  \eqref{eq:approxfn} with
  \begin{equation*}
    \delta_{n}\coleq L\,n^{-rp/\left((\alpha+1)(p+1)\right)}
  \end{equation*}
  with an $L>0$, we have
  \begin{equation*}
    \norm[\infty]{f_{n}-\ff}
    \leq
    \left(
      \left(
        C_{r}\left(p+1\right)^{1/p}\norm[\infty]{\ff}
      \right)^{p/(p+1)}\!\!\!/L
      + c\,L^{\alpha}
    \right)
    \, n^{-\alpha r p/\left((\alpha+1)(p+1)\right)}.
  \end{equation*}
  
  For exponential discretisation $X_{n}$ and $C_{\gamma}$ as in
  Corollary~\ref{cor:lpRgamman} and $f_{n}$ defined by
  \eqref{eq:approxfn} with
  \begin{equation*}
    \delta_{n}\coleq L\,\gamma^{- pn/\left((\alpha+1)(p+1)\right)},
  \end{equation*}
  with an $L>0$, we obtain
  \begin{equation*}
    \norm[\infty]{f_{n}-\ff}\leq
    \left(
      \left(
        C_{\gamma}\left(p+1\right)^{1/p}\norm[\infty]{\ff}
      \right)^{p/(p+1)}\!\!\!/L
      + c\,L^{\alpha}
    \right)
    \, \gamma^{\alpha pn/\left((\alpha+1)(p+1)\right)}.
  \end{equation*}
\end{cor}

\begin{bem}
  If $X$ is bounded and bounds for the density $f_{X}$ and its modulus
  of continuity are known explicitly, the last result is strong enough
  to construct a perfect simulation algorithm based on von Neumann's
  rejection method. Corollary~\ref{cor:boundApproxDens} can be turned
  into such an algorithm as done in 
  \citet{Dev01} for the case of infinitely
  divisible perpetuities with approximation of densities by Fourier
  inversion, \citet*{DevFilNein00} for the case of the Quicksort limit
  distribution and \citet{DevNein02} for more general fixed-point
  equations of type \eqref{glnn1}.
\end{bem}

\newcommand{\sC}{\ensuremath{s}}
\section{Algorithm and Complexity} \label{sec:complex}

In this section, we will give an algorithm for an approximation
satisfying the assumptions in the last section for many important cases.
We assume that the distributions of $A$ and $b$ are given by Skorohod
representations, i.e.~by measurable functions
$\phi,\psi:[0,1]\!\rightarrow\!\Real$, such that
\begin{equation}
  \label{eq:defPhiPsi}
  A=\phi(U) \quad\text{and}\quad b=\psi(U),
\end{equation} $U$ being uniformly distributed on $[0,1]$.
Furthermore, we assume that $\norm[\infty]\phi\leq1$ and that both 
functions are Lipschitz continuous and can be evaluated in constant
time. Now we define the discretisation $\discrete{\cdot}$ by
\begin{align}
  \discrete{Y}:= \flr{\sC(n)\, Y}\!/\sC(n),\label{eq:defU_n}
\end{align}
where $\sC(n)$ can be either polynomial, i.e.~$\sC(n)=n^{r}$ or
exponential, $\sC(n)=\gamma^{n}$.  Defining
\begin{align*}
  \An&\coleq \phi\!\left(\discrete{U}\right)\quad\text{and}\\
  \bn&\coleq \psi\!\left(\discrete{U}\right),
\end{align*}
the conditions on $\phi$ and $\psi$ ensure that 
Corollary~\ref{cor:lpR} and \ref{cor:lpRgamman} can be applied.

We keep the distribution of $\Xn$ in an array $\mathcal{A}_{n}$, where
\begin{displaymath}
  \mathcal{A}_{n}[k]\coleq\Ws{X_{n}={k}/{\sC(n)}}
\end{displaymath} 
for \mbox{$k\in\Zs$}.
Note however, that as $A$ and $b$ are bounded, $\mathcal{A}_{n}[k]=0$ at
least for $\abs{k}>\sC(n)Q_{n}$, where $Q_{n}$ can be computed recursively as
$Q_{n}=\ceil{\norm[\infty]{A}Q_{n-1}+\norm[\infty]{b}}$ and
$Q_{0}=\ceil{\norm[\infty]{X_{0}}}\!=\ceil{\Ew X}$. 

For simplicity we assume that $s(0)=s(1)=1$ and that $s(n)\in\Nat$ for
all $n$. The core of the implementation
is the following update procedure:
 \begin{algorithmic}
   \Procedure{update}{$\mathcal{A}_{n-1},\mathcal{A}_{n}$}
     \For{$i$}{$0$}{$\sC(n)-1$}
       \For{$j$}{$-\,\sC(n-1)\, Q_{n-1}$}{$\sC(n-1)\, Q_{n-1}$}
       \State \vspace{-3ex}
       \begin{align*}
         &u\quad \gets\quad \frac{i}{\sC(n)} \\[1ex]
         &k\quad \gets\quad \left\lfloor \sC(n)\,\left(
             \phi(u)\,\frac{j}{\sC(n-1)}+\psi(u)
           \right)\right\rfloor\hspace{1em} \\[1ex] 
         &\hspace{-0.3ex}\mathcal{A}_n[k]\gets\hspace{0.5ex}
         \mathcal{A}_n[k]+\frac{1}{\sC(n)}\,\mathcal{A}_{n-1}[j]
       \end{align*}
       \EndFor
     \EndFor
   \EndProcedure
 \end{algorithmic}

Furthermore, we use a procedure
$\textsc{initialize}(\mathcal{A}_{n},n)$, which creates $\mathcal{A}_{n}$ as vector
with $2\sC(n)Q_{n}$ components with $\mathcal{A}_{n}[k]=0$ for
$-\sC(n)Q_{n}\leq k \leq\sC(n)Q_{n}$.

The whole algorithm then looks like this: 
 \begin{algorithmic}
     \State \textsc{initialize}($\mathcal{A}_{0},0$)
     \vspace{-1em}
     \begin{flalign}\label{eq:algoInit}
     \mathcal{A}_{0}\!\left[\biggr.\flr{\Bigr.\sC(0)\;\Ew X}\right]\gets1&&
     \end{flalign}
     \vspace{-1.5em}
     \For{$n$}{$1$}{$N$}
       \State \textsc{initialize}($\mathcal{A}_{n},n$)
       \State \textsc{update}($\mathcal{A}_{n-1},\mathcal{A}_{n}$)
     \EndFor
     \State\Return $\mathcal{A}_{N}$
 \end{algorithmic}

Note, that \eqref{eq:algoInit} determines that we start the approximation
with $X_{0}$ as defined in \eqref{eq:defXtn}.

The complete code for polynomial discretisation for the example in
 Section~\ref{sec:keyEx}, implemented in \cpp, can be found in
 \citet{Kn06}.

To approximate the density as in \eqref{eq:approxfn} with
$\delta_{N}=d/s(N)$ for some $d\in\Nat$, we compute a new array
$\mathcal{D}_{N}$ by setting
\begin{equation*}
  \mathcal{D}_{N}[k]
  =\frac{\sC(N)}{2d}\sum_{j=k-d+1}^{k+d}\!\!\!\!\mathcal{A}_{N}[j].
\end{equation*}

To measure the complexity of our algorithm, we estimate the number of
steps needed to approximate the distribution function and the density up
to an accuracy of $1/n$. For the case that $X$ has a bounded density
$f_{X}$ which is H\"older continuous, we give asymptotic bounds for
polynomial as well as for exponential discretisation. We assume the
general condition \eqref{eq:defPhiPsi}.

\begin{lm}\label{lm:complex}
  Assume that $X$ has a bounded density $f_{X}$, which is H\"older
  continuous with exponent $\alpha\in(0,1]$.  Using polynomial
  discretisation with exponent $r$, cf. Corollary~\ref{cor:lpR}, we can
  calculate for any $n\in\Nat$ approximations $\hat{F}, \hat{f}$ of the
  distribution function $F$ and the density $f$ of $X$ with
  \begin{displaymath}
    \norm[\infty]{\hat{F}-F}\leq\frac{1}{n},\qquad
    \norm[\infty]{\hat{f}-f}\leq\frac{1}{n}
  \end{displaymath}
  in time $T_{F}(n)$ and $T_{f}(n)$ respectively with
  \begin{displaymath}
    T_{F}(n)=O\!\left(n^{(2+2/r)(p+1)/p}\right)  \quad\text{and}\quad
    T_{f}(n)=O\!\left(n^{2\left(1+1/\alpha\right)(r+1)(p+1)/(rp)}\right).
\end{displaymath}

Using exponential discretisation with parameter $\gamma$ as in
Corollary~\ref{cor:lpRgamman}, approximation to the same accuracy
takes time
\begin{displaymath}
  T'_{F}(n)=O\!\left(n^{(p+1)/p}\log{n}\right)
  \quad\text{and}\quad
  T'_{f}(n)=O\!\left(n^{(1+1/\alpha)(p+1)/p}\log{n}\right) 
\end{displaymath}
for the distribution
function and the density of $X$ respectively.

 \begin{proof}
    In one execution of {\sc update}($\mathcal{A}_{k-1},\mathcal{A}_{k}$), the
    outer loop is executed $\sC(k)$ times. The assumptions on $A$ and
    $b$ ensure that $Q_{k}=O(k)$, so we have $O(k\,\sC(k))$ runs
    of the inner loop and the whole procedure takes time
    $O\!\left(k\;\sC(k)^{2}\right)$. Hence, for any $N\in\Nat$, finding
    $\mathcal{A}_{N}$ costs time
    \begin{equation}\label{eq:runtimeUpdateN}
      O\!\left(\sum_{k=1}^{N}k\,\sC(k)^{2}\right)
      =O\!\left(N^{2}\,\sC\!\left(N\right)^{2}\right).
    \end{equation}

    For discretisations with $\sC(n)=n^{r}$ we get a running
    time of $O(N^{2r+2})$ to find $\mathcal{A}_{N}$, and \eqref{eq:kolBoundPoly} in
    Lemma~\ref{lm:kolBound} ensures that for the corresponding distribution function
    $F_{N}$ of $X_{N}$,
    \begin{displaymath}
       \norm[\infty]{F_{N}-F}\leq CN^{-rp/(p+1)}.
    \end{displaymath}
    Setting $N=\left(Cn\right)^{(p+1)/(rp)}$ and $\hat{F}\coleq F_{N}$,
    we get an approximation of the stated accuracy in time
    \begin{displaymath}
      T_{F}(n)=O(N^{2r+2})=O(n^{(2+2/r)(p+1)/p}).
    \end{displaymath}
    For the density of $X$ we use Corollary~\ref{cor:boundApproxDens}
    and $N'=\left(C'n\right)^{(\alpha+1)(p+1)/(\alpha rp)}$
    to obtain the stated bound.
    
    When using exponential discretisation, $\sC(n)=\gamma^{n}$, we need
    time $O(N^{2}\gamma^{N})$ to find $\mathcal{A}_{N}$. Using the
    corresponding results in Lemma~\ref{lm:kolBound} and
    Corollary~\ref{cor:boundApproxDens} ensures the stated running
    times.
  \end{proof}
\end{lm}

\begin{cor}\label{cor:complex}
  Assume \eqref{eq:defPhiPsi} and that $X$ has a bounded density
  $f_{X}$, which is H\"older continuous with exponent
  $\alpha\in(0,1]$. Then, using exponential discretisation as in
  Corollary~\ref{cor:lpRgamman}, approximation to an accuracy of  $1/n$
  takes time $O(n^{1+\epsilon})$ for the distribution function and time
  $O(n^{1+1/\alpha+\epsilon})$ for the density of $X$ for all
  $\epsilon>0$.
  \begin{proof}
    Note that $\norm[\infty]{\phi}\leq1$ and $\norm[p]{A}<1$ for some
    $p\geq1$ implies that $\norm[p]{A}<1$ for all $p\geq1$. Thus, in
    Lemma~\ref{lm:complex}, $p$ can be chosen arbitrarily large.
  \end{proof}
\end{cor}

\newcommand{\PX}[1][(x)]{\Pb_{X}#1}
\newcommand{\intgf}[3]{\int_{#1}^{#2}\!\!g(x,#3)\f dx}
\newcommand{\intp}[2][t]{\int_{p_{#1}}^{#2}\!\!\!g(x,#2)\f dx}
\newcommand{\intf}[3]{\int_{#1}^{#2}\!\!\!#3\f dx}
\newcommand{\fmax}[1][\gamma_{t}]{f_{\max}(#1)}
\newcommand{\maxf}{\norm[\infty]{\ff}}

\section{A simple class of perpetuities} \label{sec:ax1}

In order to make the bounds of Section~\ref{sec:convRates} explicit in
applications, we need to bound the absolute value and modulus of
continuity of the density of the fixed-point. For a simple class of
fixed-point equations, we give universal bounds in this section. For
more complicated cases, bounds have to be derived individually, which we
work out for one example in Section~\ref{sec:keyEx}.  

For fixed-point equations of the form
\begin{equation}
  \label{eq:ax+1}
  X\eqd AX+1 \qquad \text{with }A\geq0,
\end{equation}
where $A$ and $X$ are independent, we can bound the
density and modulus of continuity of $X$ using the
corresponding values of $A$.

\begin{lm}
  \label{lm:densAX+1}
  Let $X$ satisfy fixed-point equation \eqref{eq:ax+1} and $A$ have a
  density $f_{A}$.  Then $X$ has a density $f_{X}$ satisfying
  \begin{equation}\label{eq:inteqfX_AX+1}
    f_{X}(u)=\int_{1}^\infty
      \frac{1}{x}\;f_{A}\!\left(\!\frac{u-1}{x}\!\right)f_{X}(x)dx,
    \quad \text{ for }u\geq1,%
  \end{equation}%
  and $f_{X}(u)=0$ otherwise.
  \begin{proof}
    From the fixed-point equation we can see that $X\geq1$ almost
    surely. Now let $\Pb_{X}$ be the distribution of $X$. Conditioning
    on $X$, we get for any Borel set $B$:
    \begin{align*}
      \Ws{X\in B}
      &=\int_{1}^{\infty}\!\!\Ws{Ax+1\in B}d\Pb_{X}(x)\\
      &=\int_{1}^{\infty}\!\!\!\int_{B}f_{xA+1}(u)du\,d\Pb_{X}(x)\\
      &=\int_{1}^{\infty}\!\!\!\int_{B}
          \frac1x\; f_{A}\!\left(\!\frac{u-1}{x}\!\right)du\,d\Pb_{X}(x)\\
      &=\int_{B}\int_{1}^{\infty}
          \frac1x\; f_{A}\!\left(\!\frac{u-1}{x}\!\right)d\Pb_{X}(x)\;du,
    \end{align*}
    where we can use Fubini's theorem in the last step, because the
    integrand is product measurable. 
    The claim follows, as this is just the definition of a Lebesgue
    density.
  \end{proof}
\end{lm}

\begin{cor}
  Let $A$ have a bounded density $f_{A}$. Then $X$ has a density
  $f_{X}$ satisfying
  \begin{equation*}
    \norm[\infty]{f_{X}}\leq\norm[\infty]{f_{A}}.
  \end{equation*}
  \begin{proof}
    Using Lemma~\ref{lm:densAX+1} we get
    \begin{equation*}
      \norm[\infty]{f_{X}}\leq \norm[\infty]{f_{A}}\,\Erw{\frac1X},
    \end{equation*}
    but $X\geq1$ implies $\Erw{1/X}\leq1$, so the claim follows.
  \end{proof}
\end{cor}

\begin{cor}\label{cor:modcontXfromContinuousA}
  Let $A$ have a density $f_{A}$, and $\modcont_{f_A}$ be its modulus of
  continuity. Then the modulus of  continuity $\modcont_{\ff}$ of
  $\ff$ satisfies
  \begin{equation*}
    \modcont_{f_X}(\delta)\leq\modcont_{f_A}(\delta),
    \quad\delta>0.
  \end{equation*}
    \begin{proof}
    Using \eqref{eq:inteqfX_AX+1}, we obtain for any $u,v\in\Real$
    \begin{align}\label{eq:intboundMCAX+1}
      \abs{f_{X}(u)-f_{X}(v)\bigr.}\leq
      \int_{1}^{\infty}\frac1x f_{X}(x)\abs{
        f_{A}\left(\!\frac{u-1}x\right)-f_{A}\left(\!\frac{v-1}x\right)
      }dx.
    \end{align}
    
    But $x\geq1$ and the modulus of continuity $\modcont_{f_A}$ is
    monotonically increasing by definition, so we can bound
    \begin{equation*}
      \abs{f_{A}\left(\!\frac{u-1}x\right)-f_{A}\left(\!\frac{v-1}x\right)}
      \leq\modcont_{f_A}\left(\!\frac{\abs{u-v}}x\right)
      \leq\modcont_{f_A}(\abs{u-v}),
    \end{equation*}
    and plugging this into inequality \eqref{eq:intboundMCAX+1}, we
    obtain
    \begin{equation*}
      \abs{f_{X}(u)-f_{X}(v)\bigr.}
      \leq \Erw{\frac{1}{X}}\,\modcont_{f_A}(\abs{u-v}).
    \end{equation*}
    Now we use that $\Erw{1/X}\leq1$ and take the supremum over
    all suitable $u,v$.
  \end{proof}
\end{cor}

\newcommand{\fATilde}{\bar{f}_{A}}

This result is only useful if the density of $A$ is continuous, but
we can extend it to many practical examples, where $f_A$ has
jumps at points in a set $\mathcal{I}_{A}$. We use the jump function
of $f_{A}$, defined by
\begin{equation*}
  \deltaf{A}(s)=f_{A}(s)-\lim_{x\uparrow s}f_{A}(x),\quad s>0
\end{equation*}
and a modification of $f_{A}$ where we remove all jumps,
\begin{displaymath}
  \fATilde:=f_{A}
    -\sum_{s\in\mathcal{I}_{A}\setminus\{0\}}\deltaf{A}(s)\ind_{[s,\infty)}.
\end{displaymath}
Since $X\geq1$, we now denote by $\Delta_{f_X}$ the modulus of
continuity of the restriction of $f_X$ to $(1,\infty)$.

\begin{lm}
  Let $A$ have a bounded c\`adl\`ag density $f_{A}$.  Then, for all
  $\delta>0$,
  \begin{align*}
    \modcont_{f_X}(\delta)\leq
    \modcont_{\fATilde}(\delta)
    +\norm[\infty]{f_{X}}\!\!
    \sum_{s\in\mathcal{I}_{A}\setminus\{0\}}\!\!\!
      \frac{\abs{\deltaf{A}(s)}\delta}{s}.
  \end{align*}
    \begin{proof}
      We give the proof for the case that $f_{A}$ has only one jump, say
      in $s_{0}>0$. The general case then follows similarly.
      For $1\leq u<v$, we have
    \begin{equation*}
      \abs{f_{X}(u)-f_{X}(v)\bigr.}\leq
      \int_{1}^{\infty}\frac1x f_{X}(x)\abs{
        f_{A}\!\!\left(\!\frac{u-1}x\!\right)
        -f_{A}\!\!\left(\!\frac{v-1}x\!\right)
      }dx.
    \end{equation*}
    We define
    \begin{equation*}
      \alpha\coleq\frac{u-1}{s_{0}}\vee1,\quad 
      \beta\coleq\frac{v-1}{s_{0}}\vee1
    \end{equation*}
    and divide the range of integration into the three intervals
    $(1,\alpha],[\alpha,\beta],$ and $[\beta,\infty)$. Now, in the first and
    third interval, differences of values of $f_{A}$ and $\fATilde$
    coincide. Moreover, for $x\in[\alpha,\beta]$ we have
    \begin{align*}
      \abs{
        f_{A}\!\left(\!\frac{u-1}x\!\!\right)
        -f_{A}\!\left(\!\frac{v-1}x\!\!\right)
      }
      &\leq \abs{
        \fATilde\!\!\left(\!\frac{u-1}x\!\right)
        -\fATilde\!\!\left(\!\frac{v-1}x\!\right)
      }
      +\abs{\deltaf{A}(s_{0})\bigr.}.
    \end{align*}
    
    Putting everything together we obtain
    \begin{align*}
      \lefteqn{\abs{f_{X}(v)-f_{X}(u)\big.}\leq}\qquad\quad&\\
      &\leq
      \int_{1}^{\infty}\!\frac1x f_{X}(x)\abs{
       \fATilde\!\left(\!\frac{u-1}x\!\right)
        -\fATilde\!\left(\!\frac{v-1}x\!\right)
      }dx      
      + \int_{\alpha}^{\beta}\frac1x f_{X}(x)
        \abs{\deltaf{A}(s_{0})\big.}dx\\[1ex]
      &\leq
      \int_{1}^{\infty}\!\frac1x f_{X}(x)\abs{
       \fATilde\!\left(\!\frac{u-1}x\!\right)
        -\fATilde\!\left(\!\frac{v-1}x\!\right)
      }dx
      +\norm[\infty]{f_{X}}\frac{v-u}{s_{0}}\abs{\deltaf{A}(s_{0})\big.}.
    \end{align*}
    We now 
    bound the latter integral by $\modcont_{\fATilde}(v-u)$ as in
    Corollary~\ref{cor:modcontXfromContinuousA}, 
    and the claim follows by taking the supremum over all
    $v-u\leq\delta$.
  \end{proof}
\end{lm}

\section{Example: Number of key exchanges in Quickselect}
\label{sec:keyEx}

In this section, we apply our algorithm to the fixed-point
equation
\begin{equation}\label{eq:ux+u(1-u)}
  X\eqd UX+U(1-U),
\end{equation}
where $U$ and $X$ are independent and $U$ is uniformly distributed on
$[0,1]$. This equation appears in the analysis of the selection
algorithm Quickselect. The asymptotic distribution of the number of key
exchanges executed by Quickselect when acting on a random equiprobable
permutation of length $n$ and selecting an element of rank $k=o(n)$ can
be characterized by the above fixed-point equation, see
\citet{HwTsai02}.

\begin{figure}[!ht]
  \centering
  \includegraphics*{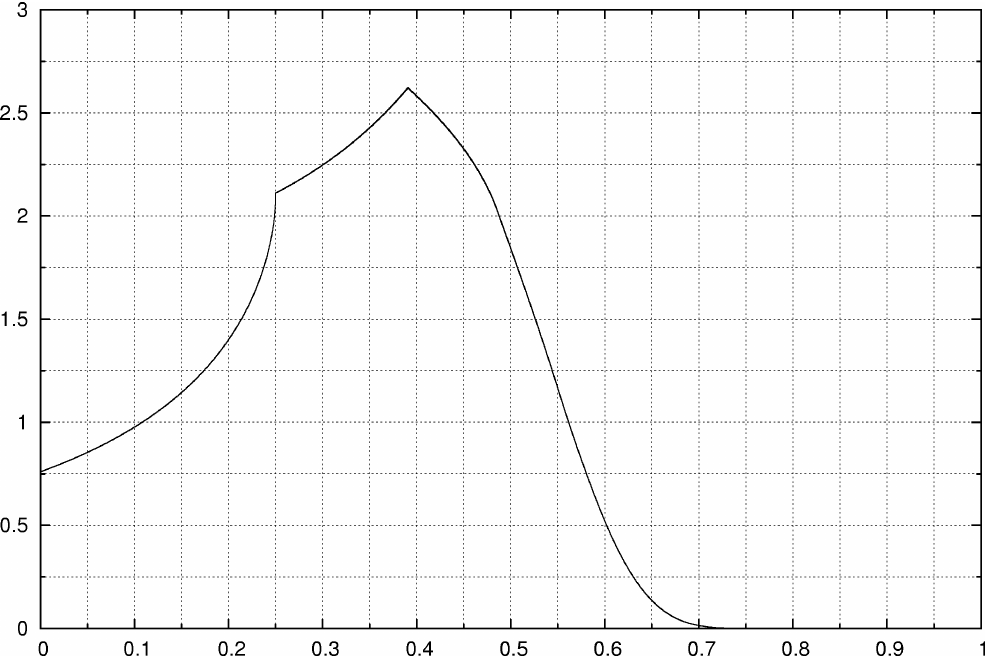}
  \caption{Histogram of approximation for $X\eqd UX+U(1-U)$.}
  \label{fig:histUX+U(1-U)}
\end{figure}

We use our algorithm to get a discrete approximation of the fixed
point. The plot of a histogram, generated with $80$ iterations of the
algorithms using for the discretisation $s(n)=n^{3}$, can be
found in Figure~\ref{fig:histUX+U(1-U)}.

In the following, we specify how the bounds in
Section~\ref{sec:convRates} can be made explicit for
this example.

\begin{lm}
  \label{lm:momBsp}
  Let $X$ be a solution of \eqref{eq:ux+u(1-u)}. Then, we have $0\leq
  X\leq1$ almost surely, and the moments are recursively given by
  $\Erw{X^{0}}=1$ and
  \begin{equation*}
    \label{eq:expXk}
    \Erw{X^{k}}
    =(k+1)!\,(k-1)!\sum_{j=0}^{k-1}\frac{\Erw{X^{j}}}{j!(2k-j+1)!},
    \quad k\geq1,
  \end{equation*}
  in particular, $\Erw{X}=1/3$.
  \begin{proof}
    
    Both claims follow directly from the fixed-point equation in
    \eqref{eq:ux+u(1-u)}, using that the solution is unique. To compute the
    moments, note that $\Erw{U^{k}(1-U)^{k-j}}$ is equal to the
    Beta function $\mathrm{B}(k+1,k-j+1)$, so we have
    \begin{align*}
      \Erw{X^{k}} 
      &=\frac{1}{1-\Erw{U^{k}}}\,
      \sum_{j=0}^{k-1}\binom{k}{j}\Erw{X^{j}}\mathrm{B}(k+1,k-j+1)\\[1ex]
      &=\frac{k+1}{k}\,\sum_{j=0}^{k-1}\frac{k!}{j!(k-j!)}
      \,\frac{k!(k-j)!}{(2k-j+1)!}\,\Erw{X^{j}}
    \end{align*}
    and the assertion follows.
  \end{proof}
\end{lm}

\begin{lm}\label{lm:Geps}
  Let $X$ be a solution of \eqref{eq:ux+u(1-u)}.
  Then, for all $\kappa\in\Nat$ and $\epsilon>0$,
  \begin{equation*}
    \Ws{X\geq1-\epsilon}
    \leq2^{(\kappa^{2}-\kappa)/4}\,\epsilon^{\kappa/2}.
  \end{equation*}
  \begin{proof}
    Using that $X$ is supported by $[0,1]$, it is easy to
    show that for all $\epsilon>0$,
    \begin{align*}
      \Ws{X\geq1-\epsilon}
      &=\Ws{UX+U(1-U)\geq1-\epsilon}\\
      &\leq\Ws{X\geq1-2\epsilon}\,\Ws{U\geq1-\sqrt{\epsilon}},
    \end{align*}
    and this inequality can be translated into
    \begin{equation}
      \label{eq:GepsFirst}
      \Ws{X\geq1-2\epsilon}\geq\frac{\Ws{X\geq1-\epsilon}}{\sqrt{\epsilon}}.
    \end{equation}

    Applying
    \eqref{eq:GepsFirst} $\kappa$ times, we get
    \begin{align*}
      1\geq 
      \Ws{X\geq1-2^{\kappa}\epsilon}
      \geq \frac{
        \Ws{X\geq1-\epsilon}}{2^{\kappa(\kappa-1)/4}
        \,\epsilon^{\kappa/2}}.
    \end{align*}
    This implies the assertion.
  \end{proof}
\end{lm}

\renewcommand{\f}[1][x]{f(#1)}
\renewcommand{\ff}{f}

\begin{lm}
  \label{lm:intEqfX}
  Let $X$ be a solution of \eqref{eq:ux+u(1-u)}. Then $X$ has a
  Lebesgue density $\ff$ satisfying $\f[t]=0$ for $t<0$ or $t>1$ and
  \begin{equation}
    \label{eq:intf}
    \f[t]=
    2\int_{p_{t}}^{t}\!\!\!g(x,t)\f dx
    + \int_{t}^{1}\!\!\!g(x,t)\f dx 
    \qquad\text{ for }t\in[0,1],
  \end{equation}
  where
  \begin{equation*}
    \label{eq:pt+gxt}
    p_{t}\coleq 2\sqrt{t}-1,
    \qquad g(x,t)\coleq\frac{1}{\sqrt{(1+x)^{2}-4t}}.
  \end{equation*}
  \begin{proof}
    Let $\PX[]$ be the distribution of $X$. Then we get for any Borel
    set $B$ by conditioning on $X$ as in the proof of
    Lemma~\ref{lm:densAX+1},
    \begin{align*}
      \Ws{X\in B}
      &=\Ws{\bigl.UX+U(1-U)\in B}\\
      &=\int_{0}^{1}\!\!\Ws{\bigl.Ux+U(1-U)\in B}d\PX\\
      &=\int_{0}^{1}\!\!\int_{B}\!\phi_{x}(t)dt\;d\PX\\
      &=\int_{B}\int_{0}^{1}\!\!\phi_{x}(t)d\PX\;dt
    \end{align*}
    where $\phi_{x}$ is a Lebesgue density of $(1+x)U-U^{2}$. The last
    step is valid by Fubini's theorem as $(x,t)\mapsto \phi_{x}(t)$ is
    product measurable, cf. \eqref{eq:defphi}. 
    
    Hence, $X$ has a Lebesgue-density $\f$ satisfying
    \begin{equation}
      \label{eq:intfphi}
      \f[t]=\int_{0}^{1}\!\!\phi_{x}(t)\f dx.
    \end{equation}
    
    To find $\phi_{x}$, we observe that $(1+x)U-U^{2}\leq (1+x)^{2}/4$
    and get
    \begin{align*}
      \lefteqn{\Ws{(1+x)U-U^{2}\leq t}=}\qquad&\\
      &=\Ws{U\leq\frac{1+x-\sqrt{(1+x)^{2}-4t}}2\quad 
        \text{or}\quad U\geq \frac{1+x+\sqrt{(1+x)^{2}-4t}}2}\\[1ex]
      &=\begin{cases}
        0 & \text{for } t<0,\\[2ex]
        \dfrac{1+x-\sqrt{(1+x)^{2}-4t}}2& \text{for }0\leq t < x,\\[2ex]
        1-\sqrt{(1+x)^{2}-4t}& \text{for }x\leq t\leq(1+x)^{2}/4,\\[2ex]
        1 & \text{otherwise.}
      \end{cases}
    \end{align*}
    
    To get a density, we differentiate with respect to $t$ and 
    rewrite as a function of $x$ yielding
    \begin{equation}
      \label{eq:defphi}
      \phi_{x}(t)=
      \begin{cases}
        \dfrac2{\sqrt{(1+x)^{2}-4t}}& \text{for }2\sqrt{t}-1< x\leq t,\\[3ex]
        \dfrac1{\sqrt{(1+x)^{2}-4t}}& \text{for }t<x\leq1,\\[3ex]
        0 & \text{otherwise.}
      \end{cases}
    \end{equation}
    Plugging this into \eqref{eq:intfphi} we get the stated integral
    equation.
  \end{proof}
\end{lm}

\begin{bem}
  \label{bem:antiderivOfg}
  The integral of $g(x,t)$ with respect to $x$ can explicitly be
  evaluated:
  \begin{equation}
    \label{eq:antiderivOfg}
    \int\!g(x,t)\,dx = \ln\left(1+x+\sqrt{(1+x)^{2}-4t}\right).
  \end{equation}
\end{bem}

\begin{bem}\label{bem:f0keyEx}
  We will see in Lemma~\ref{lm:fHoelder} that $\f$ has a version that is
  continuous on $[0,1]$. For this version we have
  \begin{equation*}
    f(0)=\Erw{\frac1{1+X}}=0.759947956\dots
  \end{equation*}
  \begin{proof}
    Using integral equation \eqref{eq:intf} we have 
    \begin{align*}
      f(0)=\int_{0}^{1}\frac{1}{1+x}\f dx,
    \end{align*}
    and by expanding the geometric series we obtain
    \begin{equation*}
      \Erw{\frac1{1+X}}=\sum_{k=0}^{\infty}(-1)^{k}\,\Erw{X^{k}},
    \end{equation*}
    which we can calculate to any accuracy using for the $k$th moments
    the formula given in Lemma~\ref{lm:momBsp}.
  \end{proof}
\end{bem}

In order to use Lemma~\ref{lm:kolBound} to bound
the deviation of our approximation, we need an
explicit bound for the density of $X$. We derive a
rather rough bound here and  see later, that we can
use the resulting bound from our approximation to
improve it.

\begin{lm}\label{lm:boundf}
  Let $\ff$ be the density of $X$ as in Lemma~\ref{lm:intEqfX}. Then
  \begin{equation*}
    \label{eq:boundf}
    \norm[\infty]{\ff}\leq 18.
  \end{equation*}
  \begin{proof}
    To get an explicit bound for $t\in[0,1]$ we simplify the integral
    equation and obtain
    \begin{align}
      \label{eq:fleq2intg}
      \f[t]&\leq2\intf{p_{t}}{1}{g(x,t)}.
    \end{align}

    We know $\f[t]$ for $t<0$, and we can bound
    $g(x,t)$, if $x$ is bounded away from $p_{t}$. Hence we split
    the integral into a left part for which we already have a bound
    for $\ff$ and a right part, in which we can bound $g$. For any
    $\gamma\in (p_{t},1]$, we have
    \begin{align}
      \label{eq:fXt}
      \f[t]
      &\leq 2 \int_{p_{t}}^{\gamma}\!\!\!g(x,t)dx
      + 2\int_{\gamma}^{1}\!\!g(x,t)\f dx,
    \end{align}
    where in the second integral, we can use that $g$ is decreasing in
    $x$ for any fixed $t$ and bound $g(x,t)\leq g(\gamma,t)$.

    For $t<1/4$, we can use that $p_{t}$ is negative, and set
    $\gamma=0$. So the first integral vanishes and only the
    second remains and we obtain 
    \begin{align}
      \f[t]&\leq2\intf{0}{1}{g(x,t)}
      \leq 2 \;g(0,t)\!\intf{0}{1}{}
      = \frac{1}{\sqrt{\tfrac14-t}}.\label{eq:fmaxleq14}
    \end{align}

    To go on, we set $\gamma=\gamma_{t}:=(p_{t}+t)/2$ and get with
    \eqref{eq:fXt}
    \begin{align*}
      \f[t]
      &\leq 2\;\mu_{t}\int_{p_{t}}^{\gamma_{t}}\!\!\!\!g(x,t)dx
      + 2\; g(\gamma_{t},t)\int_{\gamma_{t}}^{1}\!\!\!\f dx,
    \end{align*}
    where $\mu_{t}\coleq\sup\{\f[\tau]:\;\tau\in(p_{t},\gamma_{t})\}$. 

    We can calculate the first integral using the integral of
    $g$ given in \eqref{eq:antiderivOfg},
    \begin{equation}
      \label{eq:defht}
      \int_{p_{t}}^{\gamma_{t}}\!\!\!\!g(x,t)dx
      =\ln\left(
        1+\frac{(1-\sqrt{t})^{2}+(1-\sqrt{t})\sqrt{1+6\sqrt{t}+t}}
        {4\sqrt{t}}
      \right)\eqcol h(t),
    \end{equation}
    and for the second integral, 
    we obtain
    \begin{equation*}
      \int_{\gamma_{t}}^{1}\f dx 
      \leq \int_{p_{t}}^{1}\f dx
      = \Ws{X\geq1-2\!\left(1-\sqrt{t}\right)}.
    \end{equation*}

    Putting everything together we get
    \begin{align}
      \f[t]&\leq
      2\,\mu_{t}\,h(t)
      +4\;\frac{\Ws{X\geq1-2(1-\sqrt{t})}}
           {(1-\sqrt{t})\sqrt{1+6\sqrt{t}+t}}.
      \label{eq:fleqptgamma}
    \end{align}

    For $t=1/4$ we have $\gamma_{1/4}=1/8$,
    and $\mu_{1/4}\leq2\sqrt{2}$ by \eqref{eq:fmaxleq14}, so
    \begin{align}
      \f[1/4]
      \leq 4\sqrt{2}\;\ln\!\left(1+\frac{1+\sqrt{17}}{8}\right)
      +\frac{16}{\sqrt{17}}
      \leq7.\label{eq:f14leq}
    \end{align}
    From the integral equation we get for $0\leq s<t\leq1/4$ 
    \begin{align*}
      \f[t]-\f[s]
      =&\int_{0}^{1}\bigl(g(x,t)-g(x,s)\bigr)\f dx+\\
      &+\int_{0}^{s}\bigl(g(x,t)-g(x,s)\bigr)\f dx
      +\int_{s}^{t}g(x,t)\f dx\\
      >&\ 0,
    \end{align*}
    so $f$ is strictly increasing on $[0,1/4]$. Therefore, the
    bound for $t=1/4$ extends to all
    $t\in\left[0,1/4\right]\eqcol I_{0}$.
    To go on, we recursively define $b_{0}\coleq0$ and
    \begin{displaymath}
      b_{i}\coleq\left(\frac{1+b_{i-1}}{2}\right)^{2},\quad i\geq1,
    \end{displaymath}
    and 
    \begin{equation*}
      I_{2k-1}\coleq\left(b_{k},\frac{b_{k}+b_{k+1}}{2}\right],\quad
      I_{2k}\coleq\left(\frac{b_{k}+b_{k+1}}{2},b_{k+1}\right],\quad k\geq1.
    \end{equation*}
    For each interval $I_{n}$ we find a corresponding bound $M_{n}$
    for $\ff$, using that $p_{b_{i}}=b_{i-1}$ and therefore 
    $(p_{t},\gamma_{t})\subset I_{n-1}\cup I_{n-2}\text{ for }t\in I_{n}$.

    Furthermore we get for $1/4\leq t\leq1$ by differentiating the
    function $h$ defined in \eqref{eq:defht}
    \begin{align*}
      h'(t)&= c_{t}\, \Biggl(
      \frac{d}{dt}\,\frac{(1-\sqrt{t})^{2}}{4\sqrt{t}}
      +\frac{d}{dt}\,\frac{(1-\sqrt{t})\sqrt{1+6\sqrt{t}+t}}{4\sqrt{t}}
      \Biggr) ,
    \end{align*}
    where $c_{t}\geq1$. But the first summand is negative and for the
    second observe that
    \begin{align*}
        \frac{d}{dt}\,
          (1-\sqrt{t})\sqrt{1+6\sqrt{t}+t}
      &=\frac{
        \left(1-\sqrt{t}\right)\left(3+\sqrt{t}\right)
        -\left(1+6\sqrt{t}+t\right)
      }{2\,\sqrt{t}\sqrt{1+6\sqrt{t}+t}}\\[1ex]
      &=\frac{1-4\sqrt{t}-t}{\sqrt{t}\sqrt{1+6\sqrt{t}+t}}\\
      &<0,
    \end{align*}
    hence $h(t)$ is decreasing. 

    The second summand in \eqref{eq:fleqptgamma} can be bounded using
    Lemma~\ref{lm:Geps} with $\kappa=2$ yielding
    \begin{equation}
      4\;\frac{\Ws{X\geq1-2(1-\sqrt{t})}}{(1-\sqrt{t})   
        \sqrt{1+6\sqrt{t}+t}}
      \leq 4\;\frac{\Ws{X\geq1-2\,(1-\sqrt{t})}}{2\,(1-\sqrt{t})}
      \leq 4\sqrt{2}.
    \end{equation}

    So for $t\in I_{n}=(\alpha_{n},\beta_{n}]$ we have
    \begin{equation}
      \label{eq:I1}
      f(t)\leq M_{n}:=
      \left\lceil 
        2\,h(\alpha_{n})\,\max\{M_{n-1},M_{n-2}\}+4\sqrt{2}
      \right\rceil.
    \end{equation}
    Evaluating this we obtain
    \begin{align*}
      M_{0}=7,\ 
      M_{1}=13,\ 
      M_{2}=17,\  
      M_{3}=18,\ 
      M_{4}=17.
    \end{align*}
    But for $t>b_{3}$ we have $h(t)<2/7$ so the sequence
    $\left(M_{n}\right)_{n\geq0}$ is decreasing for $n\geq4$.
  \end{proof}
\end{lm}

\begin{lm}\label{lm:fHoelder}
  Let $\ff$ be the density of $X$ as in Lemma~\ref{lm:intEqfX}. Then
  $\ff$ is H\"older continuous on $[0,1]$ with H\"older exponent
  $1/2$:
  \begin{equation}
    \label{eq:modcont}
    \abs{\ff(t)-\ff(s)\bigr.}\leq 9\maxf\sqrt{t-s},
    \qquad\text{ for }0\leq s<t\leq1.
  \end{equation}
  \begin{proof}
    Using the integral equation given in Lemma~\ref{lm:intEqfX}, we
    have
    \begin{align}
      \abs{\f[t]-\f[s]\bigr.}
      &\leq 2\abs{\int_{p_{t}}^{t}\!\!g(x,t)\f dx 
        - \int_{p_{s}}^{s}\!\!g(x,s)\f dx\,} +\nonumber \\[1ex]
      &\qquad\qquad + \abs{\int_{t}^{1}\!\!g(x,t)\f 
        -\int_{s}^{1}\!\!g(x,s)\f dx}.
       \label{eq:modCont}
    \end{align}

    With explicit calculations we find
    \begin{displaymath}
      \abs{\intgf{p_{t}}{t}{t} - \intgf{p_{s}}{s}{s}}\leq 4\maxf\sqrt{t-s}
    \end{displaymath}
    and
    \begin{displaymath}
       \abs{\intgf{t}1t -\intgf{s}1s}
       \leq \maxf\sqrt{t-s}.
    \end{displaymath}
    For details see \citet{Kn06}.
  \end{proof}
\end{lm}

\begin{bem}\label{bem:notLipschitz}
  The latter lemma cannot be substantially improved, as in $t=1/4$, the
  density $\f[t]$ is not H\"older continuous with H\"older exponent
  $1/2+\epsilon$ for any $\epsilon>0$, see \citet{Kn06}.
\end{bem}

\section{Explicit error bounds for $\mathbf{X\eqd UX+U(1-U)}$}

We can now combine the bounds for the density and
its modulus of continuity with
Lemma~\ref{lm:kolBound} and Lemma~\ref{lm:difffnf}
to bound the deviation of an approximation  from
the solution of the fixed-point equation.

To approximate the density $\ff$ we  set
\begin{equation*}
  f_{n}(x)\coleq
  \begin{cases}
    \f[0] & \text{for } 0\leq x\leq\delta_{n},\\[1ex]
    \dfrac{F_{n}(x+\delta_{n})-F_{n}(x-\delta_{n})}{2\delta_{n}}
    &\text{for } \delta_{n}<x\leq1,\\
        0 & \text{otherwise,}
  \end{cases}
\end{equation*}
where $\f[0]$ is given in Remark~\ref{bem:f0keyEx} and $F_{n}$ denotes
the distribution function of $X_{n}$.

For the values used for the plot in Figure~\ref{fig:histUX+U(1-U)},
i.e.~ $s(n)=n^{3}$ and $N=80$, we can apply Corollary \ref{cor:lpR} and
obtain:
\begin{cor}
  We have $\kol(\Xn[80],X)\leq1.162\cdot 10^{-4}$, and
  \mbox{$\norm[\infty]{f_{80}-\ff}\!\leq0.931$}. Furthermore, we can
  improve the bound of Lemma~\ref{lm:boundf} and bound
  $\norm[\infty]{\ff}\!\leq3.561$.
  \begin{proof}
    We have $C_{A}=C_{b}=C_{X}=1$, hence combining
    Lemma~\ref{lm:boundf} and Lemma~\ref{lm:kolBound}, we obtain
    \begin{align*}
      \kol(\Xn,X)&\leq \left(\!  \biggl(
        \ximax^{n}\norm[p]{X}
        +\left(2+\norm[p]{X}\right)\sum_{i=0}^{n-1}
        \frac{\ximax^{i}}{\left(n-i\right)^{r}} \biggr) \,
        \left(p+1\right)^{1/p}\, \norm[\infty]{\ff}
      \right)^{p/(p+1)}\hspace{-6.5ex}.
    \end{align*}
    The moments of $X$ can be computed using
    Lemma~\ref{lm:momBsp} and we set
    $[U]_{n}\coleq\flr{n^{3}U}/n^{3}$, hence
    \begin{equation*}
      \ximax = \norm[p]{U}=\left(\dfrac1{p+1}\right)^{1/p}.
    \end{equation*}
    
    Optimizing over $p$ for $n=80$, $r=3$, and
    $\norm[\infty]{\ff}\leq18$ yields
    \begin{equation}\label{eq:boundkol}
      \kol(\Xn[80],X)\leq5.1842\cdot10^{-4}
    \end{equation}
    for $p=12$.

    Using for $\f[0]$ the value given in Remark~\ref{bem:f0keyEx}, we
    obtain for the density
        \begin{align*}
      \norm[\infty]{f_{n}-\ff} &\leq \frac1{\delta_{n}}\kol(\Xn,X) +
      9\maxf\sqrt{\delta_{n}},
    \end{align*}
    and optimizing over $\delta_{n}$, using for the Kolmogorov metric
    the bound in \eqref{eq:boundkol}, yields
    \begin{equation*}
      \norm[\infty]{f_{80}-\ff}\leq4.512
    \end{equation*}
    for $\delta_{80}=3.44\cdot10^{-4}$ (averaging $352$ values).
    
    We can now use this to improve our bound for $\norm[\infty]{\ff}$:
    Reading off the maximal value of our approximation %
    ($\norm[\infty]{f_{80}}\!\leq2.630$), we can now bound
    \begin{equation*}
      \norm[\infty]{\ff}
      \leq\norm[\infty]{f_{80}}+\norm[\infty]{f_{80}-\ff}
      \leq7.142,
    \end{equation*}
    and this in turn enables us to improve our bounds for the
    approximation, leading to
    \mbox{$\kol(X_{80},X)\leq2.2085\cdot10^{-4}$} and
    $\norm[\infty]{f_{80}-\ff}\leq1.8331$ for $\delta_{80}=3.6\cdot10^{-4}$.
    Repeating this strategy a few times, we get the stated values for
    $p=13$ and $\delta_{80}=3.7\cdot10^{-4}$ (averaging $378$ values).
  \end{proof}
\end{cor}

\begin{bem}
  Using the realistic (but yet unproven) bound of $\maxf\leq2.7$ would
  give $\kol(\Xn[80],X)\leq8.9809\cdot10^{-5}$ ($p=13$) and
  $\norm[\infty]{f_{80}-\ff}\leq0.7101$. Hence, our approach works
  well for the distribution function. However, we cannot show
  strong error bounds for the approximation of densities
  with our arguments. 
\end{bem}

However, in the next section we see that for another example the
algorithm approximates the densities much better than the error bounds
indicate.

In Table~\ref{tab:err}, the resulting error bounds for several
possible discretisations with similar running time can be found.

\begin{table}[!htb]
  \centering
  \begin{tabular}{|c|c||c|c||c|}
    \hline
    Discret.&$N$ &$\kol(\Xn[N],X)$&opt. $p$ &$\sC(N)$\\
    \hline\hline
    $n$& 22000& 0.00178& 14& 22000\\
    \hline
    $n^{2}$& 430 & 0.00025  & 16& 184900\\
    \hline
    $n^{3}$& 80 & 0.00012 & 13 & 512000\\
    \hline
    $n^{4}$& 30 & 0.00050 & 3 & 810000\\
    \hline\hline
    $1.5^{n}$ & 35 & 0.00070 & 3 & 1456110\\
    \hline
    $1.7^{n}$ & 27 & 0.00187 & 2 & 1667712\\
    \hline
  \end{tabular}
  \caption{Bounds for $\kol(\Xn,X)$ for comparable total running times
    (about 20h on a laptop computer each). 
    The discretisations are  according to Corollaries \ref{cor:lpR} and 
    \ref{cor:lpRgamman}. By $s(N)$ the number of atoms of the discrete 
    approximation is denoted, cf.~Section~\ref{sec:complex}.}
  \label{tab:err}
\end{table}

\section{An experimental view on error bounds}
\label{sec:IntSplit}

We now apply our algorithm to another fixed-point equation for which the
solution is explicitly known. We can then compare the approximation of our
algorithm with the true density and distribution function and evaluate
the actual error to get an idea of the quality of the error bounds
proven in Section~\ref{sec:convRates}. Further examples can be found in
\citet{Kn06}. It appears that the error bounds in
Section~\ref{sec:convRates} are rather loose and that the
approximation is much better than indicated by our bounds.

In the analysis of certain random interval splitting procedures the
following fixed-point equation characterizes the distribution of a point
to which a random sequence of nested intervals shrinks:
\begin{equation*}
  X\eqd \frac{1+U}2\,X + G\frac{1-U}2, 
\end{equation*}
where $G$, $U$, and $X$ are independent, $G$ is $\DBern{1/2}$ distributed
and $U$ is uniformly distributed on $[0,1]$, see \citet*{ChenGooZam84},
\citet*{ChenLinZam81}, \citet*{DevLetSes86}, and \citet{Nein01} for
details of the interval splitting context.

To approximate the fixed-point, we use a symmetric discretisation for
$(A,b)$ instead of \eqref{eq:defU_n}, setting
\begin{equation}
  \discrete{U}\coleq(2\flr{\sC(n) U}+1)/2\sC(n)
\end{equation}
and $\sC(n)=n^{3}$.

To compute the bounds as given in Section~\ref{sec:convRates}, we
can set $C_{A}=C_{b}=1/4$, $\ximax=\norm[p]{A}$, and 
$A$ is uniformly distributed on $[1/2,1]$, so
\begin{align*}
    \norm[p]{A}^{p}=\frac{2^{p+1}-1}{2^{p}\,(p+1)}
    \quad \text{ for } p\in\Nat.
\end{align*}
It is known that $X$ is $\Dbeta{2}{2}$ distributed, so we have
the moments:
\begin{align*}
  \norm[p]{X}^{p}=\prod_{s=0}^{p-1}\frac{2+s}{4+s},
  \qquad p\in\Nat.
\end{align*}
Furthermore, $X$ has the density $\f=6\,x(1-x)$, so
$\norm[\infty]{\ff}\!=1.5$.  We can now use Lemma~\ref{lm:kolBound}
and Corollary~\ref{cor:lpR} to obtain \vspace{-1ex}
\begin{align*}
  \kol(\Xn[N],X) &\leq 
  \left( 
    1.5\,\left(p+1\right)^{1/p}
    \biggl(
      \norm[p]{A}^{N}\norm[p]{X}
      +\frac{5+\norm[p]{X}}{4}
        \sum_{i=0}^{N-1}\frac{\norm[p]{A}^{i}}{\left(N-i\right)^{3}}
    \biggr)
  \right)^{\tfrac{p}{p+1}}.
\end{align*}
For $N=50$ we minimize over $p$ and get
$p_{\min}=5$ and
\begin{equation}\label{eq:kolIntr1bound}
  \kol(\Xn[50],X)\leq0.001043.
\end{equation}

As we know the limit distribution, we can read off the true error from
the output of our simulation and find
\begin{displaymath}
  \kol(\Xn[50],X)\approx0.000012.
\end{displaymath}

It is quite exactly of the order expected for a discretisation
of step size $1/n^{3}$. Note that when approximating a differentiable
function by a step function, step size and derivative impose an
unavoidable error. Comparing our approximation to a direct
discretisation by a step function of the same step size, the deviation
is at most $1.5\cdot 10^{-8}$.

Now we look at the density. 
The modulus of continuity of the density
of the $\Dbeta{2}{2}$ distribution can be bounded by
$\modcont_{f}(\epsilon)\leq6\,\epsilon$ for all positive $\epsilon$.
So for the function $f_{N}$, which we get by averaging over
$2\delta_{N}$ as in \eqref{eq:approxfn}, we get with
Lemma~\ref{lm:difffnf}
\begin{equation*}
    \norm[\infty]{f_{N}-\ff\big.}
  \leq \frac1{\delta_{N}}\kol(\Xn[N],X) + 6\, \delta_{N}.
\end{equation*}
We evaluate for $N=50$, use the bound in
\eqref{eq:kolIntr1bound}, and minimizing over $\delta_{50}$ we obtain
\begin{equation*}
  \norm[\infty]{f_{50}-\ff\big.}\leq0.1583
\end{equation*}
for $\delta_{50}=0.01318$, so we take the average over $3\,296$ values.

Reading off the true error from the simulation we obtain
\begin{displaymath}
  \norm[\infty]{(f_{n}-f)\ind_{[0.015;0.985]}}\approx0.0003
\end{displaymath}
and $\abs{f_{n}(x)-f(x)}\leq0.02$ for $x<0.015$ or
$x>0.985$. The larger errors at the boundary are caused by the
averaging procedure used to obtain $f_{n}$. 

{\bf Acknowledgements:} We thank the referee for careful reading,
pointing out some inacurracies and helping improve the presentation of
the paper.


\end{document}